\documentclass{amsart}
\usepackage{amssymb}

\newtheorem{lemma}{Lemma}[section]
\newtheorem{thm}[lemma]{Theorem}
\newtheorem{prop}[lemma]{Proposition}
\newtheorem{cor}[lemma]{Corollary}
\newtheorem{exa}[lemma]{Example}

\newcommand{\om}[1]{\Omega^{#1}}
\newcommand{\hm}[1]{\mathrm{Hom}_{\Lambda}(#1)}
\newcommand{\Ext}{\mathrm{Ext}}
\newcommand{\Ex}[2]{\Ext_{\Lambda}^{#1} (#2)}
\newcommand{\bfr}{{\mathbf r}}
\newcommand{\mfB}{\mathfrak{B}}
\newcommand{\dmn}{{\mathrm{dim}}}
\newcommand{\arr}[2]{\begin{array}{#1}#2\end{array}}

\title[Returning Arrows]{Returning Arrows for Self-injective Algebras and Artin-Schelter Regular Algebras}
\author[Guo, Yin and Zhu]{Jin Yun Guo,
Ying Yin and Can Zhu}
\email{Guo:gjy@hunnu.edu.cn\\
Yin:yinying001@hotmail.com\\ Zhu:czhu@usst.edu.cn}
\address{Jin Yun Guo and Ying Yin\\ Department of Mathematics, Key Laboratory of High Performance Computing and Stochastic Information Processing(Ministry of Education of China), Hunan Normal  University, Changsha, CHINA}
\address{Can Zhu\\College of Science\\ Shanghai University for Science and Technology\\ Shanghai, CHINA}
\thanks{This work is partly supported by
Natural Science Foundation of China \#10971172}
\dedicatory{}



\date{}

\begin{document}

\begin{abstract}
In this paper, we discuss returning arrows with respect to the Nakayama translation appearing in the quivers of some important algebras when we construct extensions.
When constructing twisted trivial extensions for a graded self-injective algebra, we show that the returning arrows appear in the quiver, that the complexity increases by $1$ in Koszul cases, and the representation dimension also increases by $1$ under certain additional conditions.
By applying Koszul duality, for each Koszul Artin-Schelter regular algebra of global dimension $l$ and Gelfand-Kirilov dimension $c$, we construct a family of Koszul Artin-Schelter regular algebras of global dimension $l+1$ and Gelfand-Kirilov dimension $c+1$, among them one is central extension and one is $l+1$-Calabi-Yau.
\end{abstract}

\maketitle

\section{Introduction}
Geiss, Leclerc and Schr\"oer have shown that for certain rigid module $I_Q$ of the preprojective algebra associated to a Dynkin quiver $Q$, the  quiver of the endomorphism algebra $\mathrm{End}(I_Q)^{op}$ is obtained from the Auslander-Reiten quiver of a path algebra $kQ$ by inserting an extra arrow $x \longrightarrow \tau x$ for each non-projective vertex $x$ along the Auslander-Reiten translation $\tau$ \cite{gls}.

Similar phenomenon also appears in Iyama's construction of $n$-complete algebras. By construction the cone of a $n$-complete algebra $\Lambda$, he gets a $n+1$-complete algebra $\Lambda'$ whose quiver is certain truncations of copies of the quivers of $\Lambda$ connected by new arrows \cite{iy}.
For certain absolute $n$-complete algebras, we show in \cite{g10} that these quivers can be embedded into the McKay quivers of some Abelian  groups and the connecting arrows can be explained as covering of the arrows of the type $x \longrightarrow \tau x$ along the Nakayama translation  $\tau$.

In fact, we have observed in \cite{gcv} this phenomenon for the McKay quivers when embedding a  finite subgroup $G$ of $\mathrm{GL}(k^m)$ into the special linear group $\mathrm{SL}(k^{m+1})$ in a some standard way, for an algebraically closed field $k$ of characteristic $0$.
The new McKay quiver is obtained from the old one by inserting arrows from $x$ to $\tau x$ for each vertex $x$, here $\tau$ is the Nakayama translation of the old McKay quiver(see \cite{gum}).

We call such new arrow from $x$ to $\tau x$ an {\em returning arrow}, since there are usually essential paths going from the vertex $\tau x$ to the vertex $x$ in the original quiver.
The appearance of the returning arrows is usually accompanied by an increase of certain numerical invariants related to the dimensions.

In this paper we discuss the returning arrows appearing when we construct twisted trivial extension of a graded self-injective algebra and when we construct generalized twisted central extension(see Section 4), that is, the extension induced by the twisted trivial extension of its Yoneda algebra, of a Koszul Artin-Schelter regular algebra.

Self-injective algebras are important and well studied in representation theory of algebras, trivial extension relates to each finite dimensional algebra a self-injective (more precisely, a symmetric) one (see, for example, \cite{ars}).
Artin-Schelter regular algebra is viewed as the homogeneous coordinate ring of a noncommutative analogue of projective space, and is studied extensively \cite{AS, ATV1, ATV2, Lev, LPWZ, SSt}.
Central extension, and more general, normal extension are two important construction for obtaining Artin-Schelter regular algebras with higher dimensions \cite{Lev, lsv}.
In Koszul cases, it is known that self-injective algebra and Artin-Schelter regular algebra are dual one another \cite{sm, m2}.

We prove in this paper a theorem on returning arrows for graded self-injective algebra $\Lambda$, that the quiver of a twisted trivial extension of $\Lambda$ is obtained from the quiver of $\Lambda$ by inserting an arrow from $x$ to $\tau x$ for each vertex $x$ of the quiver of $\Lambda$, here $\tau$ is the Nakayama translation of the quiver of $\Lambda$.
We also describes the relations for  twisted trivial extensions.

We discuss two important invariants, the complexity and the representation dimension, for the twisted trivial extensions of a Koszul self-injective algebra.
We show that for a Koszul self-injective algebra, its twisted trivial extension is Koszul and prove that if the complexity $\Lambda$ is $c$, the complexity of its twisted trivial extension  is $c+1$.
We also get lower bound for the representation dimension of the twisted trivial extension under {\bf Fg} assumption for both algebras, using Bergh's theorem \cite{ber}.
So under some stronger condition, we also shows that if the representation dimension of $\Lambda$ is $m$, the representation dimensions of a twisted trivial extension is $m+1$.

We also present a Koszul Artin-Schelter regular algebra version of our theorem on returning arrows by applying Koszul duality.
For each Koszul Artin-Schelter regular algebra $\Gamma$ of global dimension $l$ and Gelfand-Kirilov dimension $c$, we obtain a family $\{\Gamma[\sigma] | \sigma \in \mathrm{Aut}\, E(\Gamma)\}$ of Koszul Artin-Schelter regular algebra $\Gamma$ of global dimension $l+1$ and Gelfand-Kirilov dimension $c+1$.
We also single out a central extension and a Calabi-Yau algebra explicitly in this family.

The theorem on the returning arrows for graded self-injective algebras is proven in Section 2, the relations for the twisted trivial extensions are also discussed.
In in Section 3, we discuss two numerical invariants,  complexity and representation dimension, of twisted trivial extensions of Koszul self-injective algebras.
In Section 4, we apply Koszul duality and get a version of returning arrow theorem for Koszul Artin-Schelter regular algebras, we also discuss how to get a central extension and a  Calabi-Yau algebra.

\section{Returning Arrows for the Trivial Extensions}

Let $k$ be a field.
A $k$ algebra $\Gamma$ is said to be  a {\em graded algebra} if $\Gamma = \bigoplus_{t=0}^{\infty} \Gamma_t$ as vector spaces and satisfies the following conditions:
(1) $\Gamma$ is generated by $\Gamma_0$ and $\Gamma_1$,
(2) $\Gamma_0 = \bigoplus_{i=1}^n k e_i$ such that $1 = e_1 + \cdots + e_n$ is a decomposition of 1 as a sum of  orthogonal primitive idempotents.
So $\Gamma_0$ is a direct product of $n$-copies of $k$, and $\Gamma_t\Gamma_s = \Gamma_{t+s}$ for all $s,t \ge 0$.
Modules are assume to be left modules throughout this paper, when not specialized.

Let $\Lambda = \Lambda_0+ \Lambda_1+\cdots+ \Lambda_l$ be a finitely dimensional graded self-injective algebra.
Let $\bfr = \Lambda_1+ \Lambda_2+\cdots+ \Lambda_l$ be its radical.
Recall that a bound quiver $Q = (Q_0, Q_1, \rho)$ is a quiver with a set of relations, here $Q_0$ denotes the set of vertex, $Q_1$ denotes the set of arrows and $\rho$ denotes a set of relations.
Thus $\rho$ is a set of linear combinations of paths of length larger or equal to $2$.
Write $k(Q)= kQ/(\rho)$, the quotient of the path algebra $kQ$ by the ideal $(\rho)$ generated by $\rho$, and call it the {\em algebra given the bound quiver $Q$}.
A path in $Q$ is called a bound path if its image in $k(Q)$ is nonzero.
Clearly, $n=|Q_0|$ is the number of vertices of $Q$.
Assume that $\Lambda$ is the algebra given by a bound quiver $Q$, then  $1=\sum_{i=1}^n e_i$ is a sum of orthogonal primitive idempotents corresponding to the vertices of $Q$.
$\Lambda_0 $ is a semisimple algebra and $\Lambda_0 e_1, \ldots, \Lambda_0 e_{n}$ is a complete set of representatives of the isomorphic classes graded simple $\Lambda$ modules concentrated at degree zero.
$\Lambda_1$ is the vector space with $Q_1$ as a basis, and $\Lambda_t$ is the vector space spanned by the paths of length $t$.

A bound quiver $Q=(Q_0,Q_1,\rho)$ is called {\em homogeneous} provided that each of the paths appearing in a given linear combination of $\rho$ has the same length.
Clearly, $k(Q)$ is a graded algebra if $Q$ is a homogeneous bound quiver.

Fix an integer $l\ge 1$, a homogeneous bound quiver $Q$ is said to be {\em  stable of Loewy length} $l+1$ if it satisfies the following conditions:
\begin{enumerate}
\item A permutation $\tau$ is defined on of the vertex set of $Q$;

\item The maximal bound paths of $Q$ have the same length $l$;

\item  For each vertex $i$, there is a maximal bound path from $\tau i$ to $i$, and there is no bound path of length $l$ from $\tau i$ to $j$ for any $j \neq i$;

\item Any two maximal bound paths starting at the same vertex are linearly dependent.
\end{enumerate}
$\tau$  is called the {\em Nakayama translation} of the stable bound quiver $Q$.

It is shown in \cite{g10} that an algebra is a graded self-injective algebra of Loewy length $l+1$ if and only if it is given by a stable bound quiver of Loewy length $l+1$, and the Nakayama translation  is induced by the Nakayama functor

Recall that the trivial extension $\Lambda\ltimes M$ of an algebra $\Lambda$ by a $\Lambda$-bimodule $M$ is the algebra defined on the vector spaces $\Lambda\oplus M$ with the multiplication defined by $$(a,x)(b,y)=(ab, ay+xb)$$ for $a,b\in \Lambda$ and $x,y\in M$.
$\Lambda\ltimes D\Lambda$ is called the  trivial extension of $\Lambda$.

Let $\Lambda$ be a graded  self-injective algebra of Loewy length $l+1$, let $Q$ be its quiver with Nakayama translation $\tau$.
Choose a basis $\mathfrak{B}$ of $\Lambda$ such that $\mfB_t=\mfB \cap \Lambda_t$ is a basis of $\Lambda_t$ consisting of bound paths of length $t$.
Write $\mfB_t =\{p_{t,j}|j = 1,\cdots, r_t=\dmn_k \Lambda_t\}$.
We have that $\dmn_k \Lambda_l=\dmn_k \Lambda_0$ is the number of the vertices of $Q_{\Lambda}$.

For each vertex $i$, all the maximal bound paths ending at the vertex $i$ start at vertex $\tau i$ and they are pairwise linear dependent.
Write $p_{l,i}$ for the basic element which is a maximal bound path from $\tau i$ to $i$.
The dual basis $\mfB^* =\cup_{t} \mfB_t^* $  is a basis of $D\Lambda$. Denote by $\Lambda_t^*$ the space spanned by $\mfB^*_t$, then $D\Lambda = \Lambda_l^*+ \cdots +\Lambda_0^*$.

$D\Lambda$ is generated by $\Lambda_l^*$ as a $\Lambda$-bimodule.
Let $\tilde{\Lambda}=\Lambda\ltimes D\Lambda$ be the trivial extension of $\Lambda$.
As an algebra, $\tilde{\Lambda}$ is generated by $\Lambda_0$ and $\Lambda_1 + \Lambda^*_l$.
Regarding elements in $\Lambda_l^*$ as degree one elements.
$D\Lambda$ is a bimodule generated at degree $1$.
Let $\tilde{\Lambda}_0=\Lambda_0$, $\tilde{\Lambda}_t =\Lambda_t + \Lambda^*_{l-t} $ for $t=0,1,\cdots, l-1$ and $\tilde{\Lambda}_{l+1} = \Lambda_0^*$.
We see $$\tilde{\Lambda}=\tilde{\Lambda}_0+\tilde{\Lambda}_1+\cdots + \tilde{\Lambda}_l+ \tilde{\Lambda}_{l+1} $$ is also a graded self-injective algebra.

Let $Q$ and $\tilde{Q}$ be the Gabriel quiver of $\Lambda$ and of $\tilde{\Lambda}$, respectively.
We have the following lemma

\begin{lemma}\label{quiver}
The quiver  $\tilde{Q}$ is obtained from $Q$ by adding an arrow from $i$ to $\tau i$ for each vertex $i$ of $Q$.

$\tilde{Q}$ is a stable translation quiver with Loewy length $l+2$ and trivial Nakayama translation.
\end{lemma}
\begin{proof}
By definition, idempotents of $\Lambda$ and of $ \tilde{\Lambda}$ are all in $\Lambda_0$, hence they are the same.
Note $\bfr \tilde{\Lambda} = (\Lambda_1+\Lambda_l^*) + (\Lambda_2 + \Lambda_{l-1}^* )+ \cdots +( \Lambda_l + \Lambda_1^*) + \Lambda_0^*$.
We have $\tilde{\Lambda}/\bfr \tilde{\Lambda} \simeq \Lambda_0$ as semisimple algebras, so $Q$ and $\tilde{Q}$  have the same set of vertices.

On the other hand, $\bfr \tilde{\Lambda}/ \bfr^2 \tilde{\Lambda} \simeq \Lambda_1 \oplus \Lambda^*_l$ as $\Lambda_0$ bimodules.
The arrow set of $Q$ can be chosen as a subset of $\tilde{Q}$ since $\Lambda_1$ is a subspace of $\tilde{\Lambda}_1$, and the remaining arrows of $\tilde{Q}$ are afforded by a basis of $\Lambda^*_l$.
Note that a maximal bound paths $p_{l,i}$ starting at vertex $\tau i$ ends at $i$, so one gets that $$e_j p_{l,i}^* e_{i} = \left\{\arr{ll}{p_{l,i}^* & \mbox{if } j =\tau i\\ 0 & otherwise.}\right.$$
This shows that $p_{l,i}^*$ affords an arrow from $i$ to $\tau i$.
Since maximal bound paths in $Q$ ending at $i$ are pairwise linearly dependent, there is only one such new arrow for each vertex.

Clearly $\tilde{\Lambda}$ has Loewy length $l+2$. Since $\tilde{ \Lambda} $ is symmetric, the Nakayama functor is identity.
Hence the Nakayama translation is trivial.
\end{proof}

Regard $Q$ as a subquiver of $\tilde{Q}$, denote by $\alpha_i$ the new arrow $p_{l,i}^*$ in  $\tilde{Q}$ from $\tau i$ to $i$.
Now we discuss the relations for the trivial extension of a graded self-injective algebra.

\begin{lemma}\label{linf}
There is a map $\nu$ from $Q_1$ to $\cup_{i,j}e_j\Lambda_1 e_i$ such that $\nu(Q_1)$ is a basis of $\Lambda_1$, and there is a linear function $\eta_i$ on  $e_i\Lambda_{l} e_{\tau i}$ for each vertex $i$ satisfying the following conditions:
\begin{enumerate}
\item $\eta_i(p) \neq 0$ if and only if $p$ is nonzero element in $e_i\Lambda_{l} e_{\tau i}$;
\item For any arrow $\beta: i \longrightarrow j$, if $q$ is a bound path of length $l-1$, and $\beta q \neq 0$, then $\eta_j(\beta q) =\eta_i(q \nu(\beta )) $.
\end{enumerate}
\end{lemma}
\begin{proof}
Since $e_i\Lambda_{l} e_{\tau i}$ is one dimensional vector space, there is a linear function $$\eta_i: e_i\Lambda_{l} e_{\tau i} \longrightarrow k$$ on $ e_i\Lambda_{l} e_{\tau i}$ such that for any $z \in e_i\Lambda_{l} e_{\tau i}$, $z = \eta_i(z) p_{l,i}$.
For any arrow $\beta :i\longrightarrow j$, we have  a pair $B^{i,j}$ and $B_{i,j}$ of non-degenerated bilinear forms on $e_i\Lambda_{l-1} e_{\tau j} \times e_{\tau j}\Lambda_{1} e_{\tau i} $ and $e_j\Lambda_{1} e_{i} \times e_{i}\Lambda_{l-1} e_{\tau j} $, respectively, defined by $$B^{i,j}(x,y)=\eta_i(xy)\mbox{  \and  }B_{i,j}(y'x)=\eta_j(y'x)$$ for $x \in  e_{i}\Lambda_{l-1} e_{\tau j}, y\in e_{\tau j}\Lambda_{1} e_{\tau i}$, and $ y' \in e_j\Lambda_{1} e_{i} $.
Thus $ e_{i}\Lambda_{l-1} e_{\tau j}$ is the dual space of $ e_j\Lambda_{1} e_{ i}$, and $ e_{\tau j}\Lambda_{1} e_{\tau i}$ is a dual space of $e_{i}\Lambda_{l-1} e_{\tau j}$.

Let $\beta_1, \ldots, \beta_m$ be the arrows from $i$ to $j$, then they form a basis of $ e_j\Lambda_{1} e_{ i}$.
Let $\zeta_1, \ldots, \zeta_m$ be the dual basis of  $\beta_1, \ldots, \beta_m$ in $ e_{i}\Lambda_{l-1} e_{\tau j}$ and $\gamma_1, \ldots, \gamma_m$ be the dual basis of  $\zeta_1, \ldots, \zeta_m$ in $ e_{\tau j}\Lambda_{1} e_{\tau i}$.
The map $\nu_{i,j}: \beta_t  \longrightarrow \gamma_t$ for $t=1, \ldots, m$ extends to a bijective linear map from $ e_j\Lambda_{1} e_{ i}$ to $ e_{\tau j}\Lambda_{1} e_{\tau i}$.
It is obvious that $\eta_j(\beta q) =\eta_i(q \nu_{i,j}(\beta )) $ for all $q$ in $ e_{i}\Lambda_{l-1} e_{\tau j}$.
Now for any $\beta:i\longrightarrow j\in Q_1$, define $\nu(\beta) = \nu_{i,j} (\beta)$, we get the desired map $\nu$.
\end{proof}
This $\nu$ is, in fact, the Nakayama automorphism of $\Lambda$ induced by the Nakayama functor on $\mathrm{mod}\, \Lambda$.

Fix a linear function $\eta_i$ on  $e_i\Lambda_{l} e_{\tau i}$ for each vertex $i \in Q_0$.
Then for each vertex $i$, $p^*_{l,i} = c_i \eta_i$ for some $c_i \in k$, that is, the elements in the dual basis is just a multiple of $\eta_i$ of the given linear function.
We may assume that $c_i=1$ for all $i$ for the linear functions chosen in the proof of Lemma \ref{linf}.

By definition, for a given maximal bound paths $p_{l,j}$, and for any paths $q_t, q'_{t'}$ of length $t$ and
$t'$, respectively, we have that $$p_{l,j}^*(q_tq'_{t'} )=  \left\{\arr{ll}{a&q_tq'_{t'}=ap_{1,j}\mbox{ as bound paths}\\ 0 & \mbox{otherwise}}\right.$$
For any paths $q'_{t'}, q''_{t''}$, $q'_{t'} p^*_{1,j} q''_{t''}$ is defined by $$q'_{t'} p^*_{1,j} q''_{t''}(q_t) = p_{l,j}^*(q''_{t''}q_t q'_{t'} ),$$
for any path $q_t$.
We clearly have that $$q'_{t'} p^*_{1,i} q''_{t''} u'_{t'} p^*_{1,j} u''_{t''} =0$$ for any paths $q'_{t'}, q''_{t''}, u'_{t'}, u''_{t''} $, by the definition of the trivial extension.

\begin{lemma}\label{relations} In the trivial extension $\Lambda\ltimes D\Lambda$ we have:

For any vertex $i,j$, $$\alpha_j\alpha_i = 0. $$

For any $\beta: i \longrightarrow j \in Q_1$, $$\alpha_j \beta - \nu(\beta)\alpha_i = 0.$$
\end{lemma}
\begin{proof}
We clearly have that $\alpha_i \alpha_j = 0$ for all $i,j$, by the definition of the trivial extension, especially, $\alpha_{\tau i}\alpha_i = 0$.

On the other hand, for any arrow $\beta$ from $i$ to $j$, the linear functions $\alpha_j\beta$ and $\nu(\beta) \alpha_i$ are defined as follows:
$$\alpha_j \beta(q) = p^*_{l,j}(\beta q)\quad \mbox{ and }\quad\nu(\beta)\alpha_i(q) = p^*_{l,i}(q\nu(\beta))$$ for any path $q$.
So by Lemma \ref{linf}, we have $\alpha_j \beta(q)=\nu(\beta)\alpha_i(q)$.
Thus  $$\alpha_j \beta - \nu(\beta) \alpha_i =0$$ for any arrow $\beta$ in from $i$ to $j$ $Q$.
\end{proof}

Now we prove our first main theorem on returning arrows.

\begin{thm}\label{main}
Let $\Lambda$ be a graded self-injective algebra given by a bound quiver $Q=(Q, \rho)$.
Let $\tilde{\Lambda} = \Lambda \ltimes D\Lambda$ be its trivial extension, and let $\tilde{Q}=(\tilde{Q},\tilde{\rho})$ be a bound quiver of $\tilde{\Lambda}$.
If $Q$ is a stable bound quiver of Loewy length $l+1$ with the  Nakayama translation $\tau=\tau_{\Lambda}$. Then
\begin{enumerate}
\item $\tilde{Q}$ is a stable translation quiver of Loewy length $l+2$ with trivial Nakayama translation, and $\tilde{Q}$ is obtained from $Q$ by adding an arrow $\alpha_i: i \longrightarrow \tau i$ for each vertex $i\in Q_0$;

\item The relation can be taken as $\tilde{\rho}= \rho\cup \rho'$ with $$\rho'= \{\alpha_{\tau i} \alpha_i | i \in Q_0\}\cup \{\alpha_j\beta-\nu(\beta)\alpha_i| \beta: i\longrightarrow j \in Q_1\}.$$
\end{enumerate}\end{thm}

\begin{proof}
The first assertion follows from Lemma \ref{quiver}.

Now we prove that $ \tilde{\Lambda}$ is the algebra given by the bound quiver $(\tilde{Q}, \tilde{\rho})$ by proving that $k\tilde{Q}/(\tilde{\rho}) \simeq \Lambda\ltimes D\Lambda$.

Note that $kQ\cap(\tilde{\rho})=(\rho)$, so $\Lambda \simeq kQ/(\rho)$ is a subalgebra of $k\tilde{Q}/(\tilde{\rho})$.
Let $\phi: kQ / (\rho) \longrightarrow \Lambda$ be an isomorphism.
Then by Lemma \ref{quiver}  and \ref{relations},  we can extend $\phi$ to an epimorphism $\phi': k\tilde{Q}/(\tilde{\rho}) \longrightarrow \Lambda \ltimes D\Lambda$ by defining $\phi'(\alpha_i) = p^*_{l,i}$.
We now prove that $\phi'$ is an isomorphism.

Note that for each bound path $p_{t,s} \in \mfB$ from $j$ to $i$ of length $t$, there is a path $q_{t,s}$ from $\tau(i)$ to $j$ of length $l - t$, such that $p_{t,s}q_{t,s}$ is a maximal bound path.
So we see that $p_{t,s}q_{t,s}= a_{t,s}^{-1}p_{l,i}$ in $\Lambda$ for some $0 \neq a_{t,s} \in k$.
If $p_{t,s}^*$ is the corresponding basic elements in the dual basis $\mfB^*$, then we have that $p_{t,s}^* = a_{t,s}q_{t,s}p^*_{l,i}$.
Fix a $q_{t,s}$ for each $t,s$, the set $\mfB'' = \{q_{t,s}|p_{t,s} \in \mfB\}$ is also a basis of $\Lambda$ (see \cite{gmt}).

$kQ/(\rho)$ is a subalgebra of $k\tilde{Q}/(\tilde{\rho})$.
Let $\mfB' = \{q_{t,s}\alpha_i \in k\tilde{Q}/(\tilde{\rho})| p_{t,s}\in \mfB\}$, we prove that $\mfB \cup \mfB'$  spans $k\tilde{Q}/(\tilde{\rho})$.
Give a nonzero path $p= \gamma_r\cdots \gamma_1$ in $k\tilde{Q}/(\tilde{\rho})$.
If all the arrows are in $Q$, then $p \in \Lambda$ and it is a linear combination of elements in $\mfB$.
Assume that $\gamma_s = \alpha_i$ is a new arrow for $s>1$.
If $\gamma_{s-1} = \alpha_j$ is also a new arrow, then $p=0$ in $k\tilde{Q}/(\tilde{\rho})$ since $\alpha_i\alpha_j \in \tilde{\rho}$.
If $\gamma_{s-1} = \beta $ is an arrow from $j$ to $i$, then we have that $p= \gamma_r\cdots  \gamma_s \gamma_{s-1} \cdots  \gamma_1 = \gamma_r\cdots  \gamma_{s+1}\nu(\beta) \alpha_j\gamma_{s-2} \cdots  \gamma_1  $ in $k\tilde{Q}/(\tilde{\rho})$ since $\alpha_i \beta - \nu(\beta) \alpha_j \in \tilde{\rho}$.
Note that $\nu(\beta)$ is a linear combinations of arrows in $Q$.
Continue in this way, we eventually write $p = \gamma'_r \cdots \gamma'_2\alpha_{i'}$ for some $i'\in Q_0$, with $\gamma'_2, \cdots,\gamma'_r$ linear combinations of arrows in $Q$.
As an element in $kQ/(\rho)$, $\gamma'_r \cdots \gamma'_2$ is written as a linear combination of elements in $\mfB$.
Thus $p$ is written as a linear combination of elements in $\mfB'$.
So we see $k\tilde{Q}/(\tilde{\rho})$ is spanned by  $\mfB \cup \mfB'$ and
$$\dmn_k k\tilde{Q}/(\tilde{\rho}) \le |\mfB \cup \mfB'| = \dmn_k \tilde{\Lambda}.$$
Hence $k\tilde{Q}/(\tilde{\rho})$ is isomorphic to $\tilde{\Lambda}$.
\end{proof}

Let $\sigma$ be a graded automorphism of $\Lambda$, that is, an automorphism that preserves the degree of homogeneous element.
Let $M$ be a $\Lambda$-bimodule.
Define the twist $M^{\sigma}$ of $M$ as the bimodule with $M$ as the vector space.
The left multiplication is the same as $M$, and the right multiplication is twisted by $\sigma$, that is,  defined by $xb=x\sigma(b)$ for all $x \in M^{\sigma}$ and $b \in \Lambda$.
Define the twisted trivial extension $ \tilde{\Lambda}^{\sigma} = \Lambda \ltimes D \Lambda^{\sigma}$ to be the trivial extension of $\Lambda$ by the twisted $\Lambda$-bimodule $D \Lambda^{\sigma}$.
$\sigma$ restricts to an automorphism on $\Lambda_0$.
It is easy to see that $\sigma$ induces an permutation on the vertices of the quiver $Q$ which commutes with the Nakayama translation.
 So Lemma \ref{quiver} and (1) of Theorem \ref{main} remain true for such a twisted trivial extension.
\begin{prop}\label{twistquiver}
The quiver $ \tilde{Q}^{\sigma}$ of $ \tilde{\Lambda}^{\sigma}$ is $\tilde{Q}$.
\end{prop}
So $ \tilde{\Lambda}^{\sigma} \simeq k\tilde{Q} /( \tilde{ \rho }^{\sigma})$.
We need to twist the relation $\tilde{\rho}$ to get $ \tilde{ \rho }^{\sigma}$.
By (2) of Theorem \ref{main}, the following proposition follows easily.

\begin{prop}\label{twistR}
The relation $ \tilde{\rho}^{\sigma}$ of $ \tilde{\Lambda}^{\sigma}$ can be taken as $\tilde{\rho}^{\sigma}= \rho\cup \rho''$ with $$\rho''= \{\alpha_{\tau i} \alpha_i | i \in Q_0\}\cup \{\alpha_j\beta-\nu\sigma^{-1}(\beta)\alpha_i| \beta: i\longrightarrow j\in Q_1\}.$$
\end{prop}

Since $\Lambda$ is graded, the linear map $x \longrightarrow (-1)x$ on $\Lambda_1$ induces a graded automorphism $\varepsilon$ on $\Lambda$.
Write $\varepsilon' = \varepsilon^{l+1}$, it sends $x \longrightarrow (-1)^{l+1} x$.
\begin{exa}\label{gradeds}
{\em Note that $\tilde{\Lambda}$ is a symmetric algebra of Loewey length $l+2$.
The non-degenerate associative bilinear form $(\,\, ,\,\, ): \tilde{ \Lambda} \times \tilde{\Lambda} \longrightarrow k$ is induced by the composition of the multiplicative map with the projection to the socle.
So $(x,y)=0$ for any  $x, y \in \tilde{\Lambda}$ with the degree $l+1$ component $(xy)_{l+1} =0$.
Let $\sigma= \varepsilon'$.
For the induced non-degenerate associative bilinear form
$$(\,\, ,\,\, )_{\sigma}: \tilde{\Lambda}^{ \varepsilon'} \times \tilde{ \Lambda}^{ \varepsilon'} \longrightarrow k$$ of $\tilde{\Lambda}^{ \varepsilon'}$, we have that $(x',y')=(-1)^{st}(y',x')$ for $x'\in \tilde{\Lambda}^{ \varepsilon'}_s, y'\in \tilde{\Lambda}^{ \varepsilon'}_t$.
This shows that $\tilde{\Lambda}^{ \varepsilon'}$ is a graded symmetric algebra.

We remark that $ \tilde{ \Lambda }^{ \varepsilon'} = \tilde{ \Lambda } $ when $l$ is odd.
}\end{exa}

\begin{exa}\label{skewcentral}{\em
Let $\sigma_0= \varepsilon \nu $, then in $\tilde{ \Lambda }^{ \sigma_0 }$, the new generators $\alpha_i$ satisfies the skew commutative relations $\alpha_j\beta + \beta\alpha_i$ for any arrow $\beta:i\longrightarrow j$ in $Q$.

If $\Lambda$ is exterior algebra $\wedge V$ of an $m$-dimensional vector space, its quiver $Q$ has one vertex and $m$ loops $x_1,\ldots, x_m$. So the quiver $\tilde{Q}^{\sigma_0}$ has one more loop (the returning arrow) $x_{m+1}$ with additional relation $x_{m+1}^2$ and $x_{m+1}x_t+x_t x_{m+1}$ for $t=1,\ldots,m$. So $\tilde{\Lambda}^{\sigma_0}$ is the exterior algebra of an $m+1$--dimensional vector space.
}
\end{exa}

\section{Complexity and Representation Dimension}

Complexity and representation dimension are important numerical invariants for algebras.
The complexity is the same as  the Gelfand-Kirilov dimension of its Yoneda algebra, an Artin-Schelter regular algebra, when the algebra is Koszul self-injective \cite{g4}.

The {\em complexity $c(\{b_r\})$} of a sequence $\{b_r\}_r$ of positive numbers is  the least non-negative number $d$ such that there exists $\lambda >0$ with $c_r \le \lambda \cdot r^{d-1}$ for almost all natural numbers $r$.
Let $M$ be a $\Lambda$-module and $$P^{\bullet}(M):\qquad \ldots  \longrightarrow P^{r}(M) \longrightarrow \cdots \longrightarrow P^{0}(M) \longrightarrow M \longrightarrow 0 $$
be a minimal projective resolution of  $M$ and let $\Omega^{r}M$ be the $r$th syzygy of $M$.
Then  {\em complexity $c(M)$} of $M$ is defined as the  complexity $c(\{\dmn_k P^{r}\}_r)$ of the sequence $\{\dmn_k P^{r}\}_r$.
Clearly, $c(M)=c(\{\dmn_k \Omega^{r}M\}_r)$.
The  {\em complexity $c(\Lambda)$} of the algebra $\Lambda$ is defined as the supremum of the complexities of finite generated $\Lambda$-modules.

%By a twisted  trivial extension we always mean the twisted trivial extension $\tilde{\Lambda}^{\sigma}$ defined at the end of the last section.
We will prove the following theorem in this section.

\begin{thm}\label{complexity}
Let $\Lambda$ be a Koszul self-injective algebra and $\tilde{\Lambda} = \Lambda\ltimes D\Lambda$ be the trivial extension of $\Lambda$.
Then $$c( \tilde{\Lambda}) = c(\Lambda)  + 1.$$
\end{thm}

Let $\Gamma = \bigoplus_{t=0}^{\infty} \Gamma_t$ be a graded algebra, recall the Yodeda algebra $E(\Gamma)$ of $\Gamma$ is the vector space $\coprod_{t \ge 0} \Ext^t_{\Gamma}(\Gamma_0, \Gamma_0)$ with Yoneda product as its multiplication.
Let $\Lambda$ be a Koszul self-injective algebra of Loewy length $l+1$ and complexity $c$, then  its Yoneda algebra $\Gamma = E(\Lambda) = \coprod_{t \ge 0} \Ext^t_{\Lambda}(\Lambda/ \bfr, \Lambda/ \bfr)$ is an Artin-Shelter regular algebra (\cite{sm, m2}) of global dimension $l$ and Gelfand-Kirilov dimension  $c$ (\cite{g4}).
Both $\Lambda$ and $\Gamma$ have the same Gabriel quiver $Q$.
Their relations span a pair of dual subspaces in the subspace of the path algebra $kQ$ spanned by the paths of length $2$ (\cite{bgs,gm}).

We first consider the Kosulity of trivial extensions.

\begin{lemma}\label{klity} Let $\Lambda$ be a Koszul algebra and let $M$ be a $\Lambda$-bimodule which is projective as a left and as right $\Lambda$-modules, if $M$ is generated at degree $1$, then the trivial extension $\Lambda \ltimes M$ is Koszul.
\end{lemma}
\begin{proof}
Write $\tilde{\Lambda}^M = \Lambda \ltimes M$.
For any $\Lambda$ module $L$, the $\tilde{\Lambda}^M$ module $\tilde{ \Lambda }^M \otimes_{\Lambda} L =  L+ M \otimes_{\Lambda} L = L \oplus M \otimes_{\Lambda} L $ as $\Lambda$-module ($\oplus$ in this proof means a direct sum of $\Lambda$-modules).
Let $N$ be a finite generated left $\Lambda$-module.
Assume that
$$P^{\bullet}(N):\qquad \ldots  \longrightarrow P^{r}(N) \stackrel{f_r}{\longrightarrow} \cdots \stackrel{f_1}{\longrightarrow} P^{0}(N) \stackrel{f_0}{\longrightarrow} N \longrightarrow 0 $$
is a minimal projective resolution of $N$ as a $\Lambda$-module.
We have that $0 \longrightarrow \Omega_{\Lambda} N \longrightarrow P^0(N) \longrightarrow N \longrightarrow 0$ is the projective cover of $N$ as $\Lambda$-modules.
So $\tilde{\Lambda}^M \otimes_{\Lambda}P^{0}(N) \longrightarrow N \longrightarrow 0$ is the projective cover as $\Lambda\ltimes M$-modules, $\Omega_{\Lambda} N \oplus M \otimes P^{0}(N)$ is a $\tilde{\Lambda}^M $-submodule of $\tilde{\Lambda}^M \otimes_{\Lambda} P$, and we have a short exact sequences
$$0 \longrightarrow  \Omega_{\Lambda} N   \oplus M \otimes_{\Lambda} P^{0}(N)  \longrightarrow (\tilde{\Lambda}^M ) \otimes_{\Lambda}P^{0}(N) \longrightarrow N \longrightarrow 0, $$
and $\Omega_{\tilde{\Lambda}^M } N \simeq  \Omega_{\Lambda} N   \oplus M \otimes_{\Lambda} P^{0}(N) $.
We also have a projective cover
$$ P^1(N)\oplus (\tilde{\Lambda}^M ) \otimes_{\Lambda} M \otimes_{\Lambda} P^{0}(N) \longrightarrow  \Omega N   + M \otimes_{\Lambda} P^{0}(N)\longrightarrow 0
$$
as $\tilde{\Lambda}^M$-modules.

Assume that
$$\Omega^r_{\tilde{\Lambda}^M } N = \Omega^{r}_{\Lambda} N \oplus \sum_{j=0}^{r-1} M^{\otimes_{\Lambda} j+1}\otimes_{\Lambda} P_{r-1-j}$$
as $\Lambda$-modules, then its projective cover as $\tilde{\Lambda}^M $ module is
$$ \tilde{\Lambda}^M \otimes_{\Lambda}(\sum_{j=0}^{r} M^{\otimes_{\Lambda} j}\otimes_{\Lambda} P_{r-j}).$$
So we get a short exact sequence
$$0 \longrightarrow  \Omega^{r+1}_{\Lambda} N \oplus\sum_{j=0}^{r} M^{\otimes_{\Lambda} j+1}\otimes_{\Lambda} P_{r-j} \longrightarrow  \tilde{P}^{r+1}(N) \longrightarrow \Omega^{r}_{\tilde{\Lambda}^M } N\longrightarrow 0,$$
and $\Omega^{r}_{\tilde{\Lambda}^M } N =\Omega^{r+1}_{\Lambda} N \oplus\sum_{j=0}^{r} M^{\otimes_{\Lambda} j+1}\otimes_{\Lambda} P_{r-j}$.
Let
$$\tilde{P}^{r}(N) = \tilde{\Lambda}^M \otimes_{\Lambda}(\sum_{j=0}^{r} M^{\otimes_{\Lambda} j}\otimes_{\Lambda} P_{r-j})$$
for the $\Lambda$-module $\sum_{j=0}^{r} M^{\otimes_{\Lambda} j}\otimes_{\Lambda} P_{r-j}$.
It follows from induction that
$$\tilde{P}^{\bullet}(N):\qquad \ldots  \longrightarrow \tilde{P}^{r}(N) \stackrel{\tilde{f_r}}{\longrightarrow} \cdots \stackrel{\tilde{f_1}}{\longrightarrow} \tilde{P}^{0}(N) \stackrel{\tilde{f_0}}{\longrightarrow} N \longrightarrow 0 $$ is a projective resolution of $N$ as a $\tilde{\Lambda}^M $-module.

By our assumption $M$ is generated at degree $1$, if $P_t$ is generated at degree $t$, then $\tilde{P}^{r}(N)$ is generated at degree $r$, for $r=0,1, \ldots$.

This shows that if $N$ is Koszul as a $\Lambda$-module, then it is also  Koszul as a $\Lambda\ltimes M$-module.
Especially,  if $\Lambda$ is a Koszul algebra, so is  $\Lambda\ltimes M$.
\end{proof}

%Let $\Lambda$ be a Koszul self-injective algebra and let $\sigma$ be a graded automorphism of $\Lambda$.
As a corollary, we have the following result.
\begin{cor}\label{Koszulslf}
Let $\Lambda$ be a Koszul self-injective algebra, and let $\sigma$ be a graded automorphism of $\Lambda$. Set the dual basic elements of the maximal paths of $\Lambda$ as degree $1$ generators of $D\Lambda^{\sigma}$. Then the twisted trivial extension $\tilde{\Lambda}^{\sigma}= \Lambda \ltimes D\Lambda^{\sigma}$ is a Koszul algebra.
\end{cor}
\begin{proof}
Since $\Lambda$ is self-injective, $D\Lambda^{\sigma}$ is projective as left $\Lambda$-module and as right $\Lambda$-module.
The corollary follows immediately from Lemma \ref{klity}.
\end{proof}

Let $\Lambda$ be a graded algebra with quiver $Q= (Q_0,Q_1)$, and let $n=|Q_0|$.
The Hilbert polynomial $H(\Lambda,t)$ of $\Lambda$ is defined as $n\times n$-matrix with entries in power series ring $\mathbb Z[[t]]$:
$$H(\Lambda,t)_{i,j} = \sum_r t^r \dmn_k(e_i \Lambda_r e_j).$$
This can also be regarded as a formal power series with coefficients in the ring of integral $n\times n$  matrices.
It follows from Theorem 2.11.1 of \cite{bgs} that $$H(\Lambda,-t)H(\Gamma,t) = E,$$ where $E$ is the $n \times n$ identity matrix, if $\Lambda$ is a Koszul algebra and $\Gamma$ its Yoneda algebra.
So we write $$H(\Gamma,t)=H(\Lambda,-t)^{-1}.$$
By the Koszul duality, we can relate $H(\Gamma,t) =H(\Lambda,-t)^{-1}$ to  the minimal  projective resolution of the $\Lambda$-module $\Lambda_0$.

\begin{prop}\label{inverse}
Let $\Lambda$ be a finite dimensional Koszul algebra, and let $H(\Lambda, t)$ be its Hilbert polynomial.
Assume that $$H(\Lambda,-t)^{-1} = A_0+A_1t^1 + \cdots + A_r t^r + \cdots,$$
where $A_r$ are $n \times n$ integral matrices with nonnegative entries.
If
$$P^{\bullet}(i):\qquad \ldots  \longrightarrow P^{r}(i) \longrightarrow \cdots \longrightarrow P^{0}(i) \longrightarrow \Lambda_0 e_i \longrightarrow 0$$
is a minimal projective resolution of $\Lambda_0 e_i$ and $$ P^{r}(i) \simeq \bigoplus_{j} a_{ji}^{(r)}\Lambda e_j.$$
Then $$A^r=(a_{ij}^{(r)}), $$ for $r= 0, 1, \cdots.$
\end{prop}
\begin{proof}
Recall $\Gamma=E(\Lambda)= \bigoplus_{r \ge 0} \Ex{r}{\Lambda_0,\Lambda_0} $ is the Yoneda algebra.
And $\Lambda \simeq E(\Gamma)$.
The homogeneous component of $\Gamma$ of degree $r$ is $\Gamma_r = \Ex{r}{\Lambda_0,\Lambda_0} $ and $\Gamma =\Gamma_0+\Gamma_1 + \cdots$. Write $J =\Gamma_1 + \Gamma_2 +\cdots$, then $J^r=\Gamma_r + \Gamma_{r+1} +\cdots$, and $J_r/J_{r+1} \simeq \Gamma_r$ as $\Gamma_0= \Lambda_0$ bimodules.
Since $\Lambda_0$ is semisimple, we have $\Ex{r}{\Lambda_0,\Lambda_0} \simeq \hm{\om{r}\Lambda_0,\Lambda_0}=\hm{P^{r}(\Lambda_0),\Lambda_0}$ for each $r \ge 0$.

By \cite{gm}, $F = \coprod_{t \ge 0} \Ext^t_{\Lambda}(-, \Lambda/ \bfr)$ and $G = \coprod_{t \ge 0}\Ext^t_{\Gamma}(-, \Gamma /J )$ are duality between the categories of Koszul modules and degree $0$ morphisms in $\mathrm{mod}\, \Lambda$ and $\mathrm{mod}\, \Gamma$.
Identify the primitive idempotents of $\Lambda$ and $\Gamma$, we have that $F(\Lambda_0 e_i) =\Gamma  e_i$ \cite{gm}.
Let $\om{r} M$ be the $r$th syzygy of the Koszul module $M$, and $P^r(M)$ be the $r$th term in the minimal projective resolution of $M$.
For a Koszul $\Lambda$ module $M$, we have that $F(P^r(M)) = F(M)_r$  by \cite{gm}.
Now assume that $P^r(\Lambda e_i) = \bigoplus_j a^{(r)}_{i,j} \Lambda e_j.$
Then  $$\Gamma_r e_i =   F(P^r(\Lambda_0 e_i)) = \bigoplus_j a^{(r)}_{i,j} \Gamma_0 e_j$$ and $\dmn_k e_j \Gamma_r e_i = a^{(r)}_{ij}$.
So by the definition of Hilbert polynomial, we get $A^r=(a_{ij}^{(r)})$, for $r= 0, 1, \cdots$.
\end{proof}

For a matrix $A=(a_{ij})$, write $\|A\|   = \sum_{i,j} a_{ij}$.
Let $\Lambda$ be a Koszul self-injective algebra of Loewy length $l+1$, let $H(\Lambda, t)$ be its Hilbert polynomial.
Write $H(\Lambda, -t)^{-1} = \sum_r A_r t^r$.
Then by Proposition \ref{inverse}, all the matrices $A_r$ have nonnegative integral entries and we have the following  proposition.

\begin{prop}\label{cptymat}
The complexity of $\Lambda$ is exactly $c(\{\|A_r\|\})$.
\end{prop}

The following lemma is well known.

\begin{lemma}\label{dual}  For all $i,j \in Q_0$ and $r=0,1,\ldots,l$,
$$dim_k e_i (D\Lambda_{r}) e_j = dim_k D(e_j\Lambda_{r} e_i) %= dim_k D(e_j\Lambda_{r} e_i)
.$$
\end{lemma}

The following lemma is follows from the proof of Lemma \ref{linf} (see also Lemma 2.1 of \cite{g4}).

\begin{lemma}\label{Nakatrans}  For all $i,j \in Q_0$ and $r=0,1,\ldots,l$,
$$dim_k\, e_{ i}\Lambda_{l-r} e_j = dim_k\,e_j\Lambda_{r} e_{\tau i}.$$
\end{lemma}

Now define $P=(p_{i,j})$ to be the matrix with  $p_{i,\tau i}= 1$ for $1\leq i\leq n$, and $p_{i,j}=0$ for $j \neq \tau i$.
$P$ is a permutation matrix.
We have the following theorem on the Hilbert polynomial of the trivial extension of a Koszul self-injective algebra.

\begin{thm}\label{HilbertOfTE} Assume that $\Lambda$ is a Koszul self-injective algebra,  and let $\tilde{\Lambda}= \Lambda \ltimes D\Lambda$ be the trivial extension of $\Lambda$.
Then %there is an permutation matrix $P$ such that
$$H(\tilde{\Lambda},t) = H(\Lambda,t)(E+tP).$$
\end{thm}
\begin{proof}
We have that $\Lambda = \Lambda_0 + \Lambda_1+ \cdots+ \Lambda_{l} $ and $D\Lambda = \Lambda_l^* + \Lambda_{l-1}^*+ \cdots+ \Lambda_{0}^* $.
Thus $$\tilde{\Lambda} = \Lambda\ltimes D\Lambda = \tilde{\Lambda}_0 + \tilde{\Lambda}_1+ \cdots+ \tilde{\Lambda}_{l} +\tilde{\Lambda}_{l+1},$$
where $\tilde{\Lambda}_0=\Lambda_0$, $\tilde{\Lambda}_r =\Lambda_r + \Lambda^*_{l+1-r} $ for $r=1,2,\cdots, l$ and $\tilde{\Lambda}_{l+1} =\Lambda^*_{0} $.
So the $i,j$ entry of $H(\tilde{\Lambda},t)$ is: $$H(\tilde{\Lambda},t)_{i,j} = H(\Lambda,t)_{i,j} +t ( \sum\limits_{r=1}^{l} t^{r-1}\dmn_k(e_i\Lambda^*_{l+1-r} e_j)+t^l\dmn_k(e_i \Lambda_0^* e_j)).$$
By Lemma \ref{dual} and  \ref{Nakatrans},
$$dim_k\, e_i \Lambda^*_{r} e_j = dim_k\, D(e_j\Lambda_{r} e_i) = dim_k\, e_j\Lambda_{r} e_i =dim_k\, e_i\Lambda_{l-r} e_{\tau j}.$$
This shows that $$\arr{rcl}{H(\tilde{\Lambda},t)_{i,j} &= &H(\Lambda,t)_{i,j}+t(\sum\limits_{r=1}^{l} t^{r-1}\dmn_k \, e_{i}\Lambda_{r-1} e_{\tau j}+t^l\dmn_k\, e_i \Lambda_l e_{\tau j})\\
&=& H(\Lambda,t)_{i,j}+t(\sum\limits_{r=0}^{l} t^{r}\dmn_k\, e_{i}\Lambda_{r} e_{\tau j})\\
&=&(H(\Lambda,t)(E+tP))_{i,j}.}$$
%Clearly
%$$(H(\Lambda,t)(E+tP))_{i,j} = H(\Lambda,t)_{i,j} + t(\sum\limits_{n=0}^{l} t^{r}\dmn_k(e_i\Lambda_{r} e_{\tau j})).$$
Hence $H(\tilde{\Lambda},t)=H(\Lambda,t)(E+tP).$
\end{proof}

The next lemma follows from Proposition 2.9 of \cite{g4} (see also Theorem 2.6 of \cite{dt}).

\begin{lemma}\label{bound} Let $\Lambda$ be a Koszul self-injective algebra of complexity $d$. Then there are positive numbers $\lambda_1, \lambda_2$ such that $$\lambda_1 r^{d-1}\leq \|A_{r}\|\leq \lambda_2 r^{d-1}$$
for almost all $r$.
\end{lemma}

\begin{lemma}\label{cptyTE} If $H(\Lambda,-t)^{-1} = A_0+A_1t^1 + \cdots + A_r t^r + \cdots,$ where $A_r$ are $n \times n$ matrices with nonnegative entries.
Then complexity of $\tilde{\Lambda}$ is $$c(\tilde{\Lambda}) = c(\{\sum\limits_{s=0}^{r} \|A_s\|\}_r).$$
\end{lemma}
\begin{proof} Let $$H(\tilde{\Lambda},-t)^{-1}= \sum_{r} A_r{'}t^r = A_0^{'} + A_1^{'}t + A_2^{'}t^{2} + \cdots + A_r^{'}t^{r} + \cdots.$$
Since $(E-tP)^{-1} = ( 1 + t + tP + t^{2}P^{2} + \cdots)$, by Theorem \ref{HilbertOfTE}, we get
$$ \arr{ccl}{H(\tilde{\Lambda},-t)^{-1} &=& (E-tP)^{-1}H(\Lambda,-t)^{-1}\\ & =& ( 1 + t + tP + t^{2}P^{2} + \cdots)(A_0+A_1 t+A_2 t^{2}+\cdots+A_r t^{r} + \cdots ).}$$
This shows that $A_r^{'} = P^{r} A_0 + P^{r-1} A_1  + \cdots + A_r$.
Since $P$ is a permutation matrix, we have that $\|P^s A\| = \| A \|$ for any matrix $A$ and nonnegative integer $s$.
So $$\|A_r^{'}\| = \|A_0P^{r}\| + \|A_1P^{r-1}\| + \cdots + \|A_r\| = \|A_0\| + \|A_1\| + \cdots + \|A_r\|.$$
This proves that $$ c(\tilde{\Lambda}) = c(\{\sum\limits_{s=1}^{r} \|A_s\|\}).$$
\end{proof}

\begin{proof}[Proof of Theorem \ref{complexity}]
Since $c(\Lambda)=c(\{\|A_r\|\})=d$, so by Lemma \ref{bound}, there are  positive number $\lambda'$ and $\lambda$ and a positive integer $N$, such that $$\lambda' s^{d-1} \leq \|A_s\|\leq \lambda s^{d-1}$$ for $s> N$.
Let $\sum\limits_{s=1}^{N} \|A_s\|= h$.
So $$\sum\limits_{s=1}^{r} \|A_s\|\leq h + \lambda (\sum\limits_{u=N+1}^{r} u^{d-1}) \thicksim O(r^{d}).$$
This shows that $$c(\{\sum\limits_{r=1} \|A_r\|\})\leqq d+1.$$

On the other hand, $$\sum\limits_{s=1}^{r} \|A_s\| \ge \sum \limits_{s=1}^{N} \|A_s\| + \lambda_1 \sum\limits_{s=N+1}^{r} r^{d-1}. $$
So by Lemma \ref{bound}, there is $\lambda^{'}>0$, such that $$\sum\limits_{s=1}^{r} \|A_s\|\ge  \lambda^{'} r^{d},$$
for sufficient large $r$.
This proves that $$c(\{\sum\limits_{s=1} ^{r} \|A_s\|\}_r)\ge d+1.$$
Hence
$$c(\tilde{\Lambda}) = c(\{\sum\limits_{s=1}^{r} \|A_s\|\}_r)=d+1 =c({\Lambda})+1.$$
\end{proof}

Let $\sigma$ be a graded automorphism of $\Lambda$, note that $\sigma$ induces a permutation on $Q_0$, denote as $\sigma$ too. It is easy to see that $\sigma$ commutes with the Nakayama translation $\tau$.
Consider the twisted trivial extension $\tilde{\Lambda} = \Lambda\ltimes D\Lambda^{\sigma}$.
Then $$\tilde{\Lambda}^{\sigma} = \Lambda\ltimes D\Lambda^{\sigma} = \tilde{\Lambda}^{\sigma}_0 + \tilde{\Lambda}^{\sigma}_1+ \cdots+ \tilde{\Lambda}^{\sigma}_{l} +\tilde{\Lambda}^{\sigma}_{l+1},$$
where $\tilde{\Lambda}_0=\Lambda_0$, $\tilde{\Lambda}_r =\Lambda_r + \Lambda^{*\sigma}_{l+1-r} $ for $r=1,2,\cdots, l$ and $\tilde{\Lambda}_{l+1} =\Lambda^{*\sigma}_{0}. $
So we have that
$$\arr{rcl}{H(\tilde{\Lambda}^{\sigma},t)_{i,j}
&=&H(\Lambda,t)_{i,j} +t ( \sum\limits_{r=1}^{l} t^{r-1}\dmn_k(e_i\Lambda^{*\sigma}_{l+1-r} e_j)+t^l\dmn_k(e_i \Lambda^{*\sigma}a_0 e_j))\\
&=&H(\Lambda,t)_{i,j} +t ( \sum\limits_{r=1}^{l} t^{r-1}\dmn_k(e_i\Lambda^*_{l+1-r} \sigma(e_j))+t^l\dmn_k(e_i \Lambda_0^* \sigma(e_j)))\\
&= &H(\Lambda,t)_{i,j}+t(\sum\limits_{r=1}^{l} t^{r-1}\dmn_k \, e_{i}\Lambda_{r-1} e_{\tau \sigma j}+t^l\dmn_k\, e_i \Lambda_l e_{\tau\sigma j})\\
&=& H(\Lambda,t)_{i,j}+t(\sum\limits_{r=0}^{l} t^{r}\dmn_k\, e_{i}\Lambda_{r} e_{\tau \sigma j}).}$$

So we have the following twisted version of Theorem \ref{HilbertOfTE} and Theorem \ref{complexity}.

\begin{thm}\label{tHilbertOfTE} Let $\sigma$ be a graded automorphism of Koszul self-injective algebra $\Lambda$,  and let $\tilde{\Lambda}^{\sigma}= \Lambda \ltimes D\Lambda^{\sigma}$ be the twisted trivial extension of $\Lambda$.
Then %there is an permutation matrix $P$ such that
$$H(\tilde{\Lambda},t) = H(\Lambda,t)(E+tP^{\sigma}),$$
where $P^{\sigma}=(p_{ij})$ is the permutation matrix with $p_{i,\tau \sigma i}=1$ for $1\leq i\leq n$, and $p_{i,j}=0$ for $j \neq \tau \sigma i$.
\end{thm}

\begin{thm}\label{tcomplexity}
Let $\Lambda$ be a Koszul self-injective algebra and let  $\sigma$ be a graded automorphism of $\Lambda$.
Let $\tilde{\Lambda}^{\sigma}= \Lambda\ltimes D\Lambda^{\sigma}$ be the twisted trivial extension of $\Lambda$.
Then $$c( \tilde{\Lambda}^{\sigma}) = c(\Lambda)  + 1.$$
\end{thm}

Denote by the representation dimension of an algebra $\Lambda$.
Recall that for a semisimple algebra $\Lambda$, $\mathrm{ repdim }\, \Lambda = 1$, and for a non-semisimple  algebra $\Lambda$, $\mathrm{ repdim }\, \Lambda$ is the smallest number $d$ which can be realized as the global dimension of the endomorphism ring of a $\Lambda$-module $M$ which is both a generator and a cogenerator.

In \cite{ber}, the following assumption is induced:

{\bf Assumption (Fg).} There is a commutative Noetherian graded $k$ algebra $H = \bigoplus_{i=0}^{\infty} H_i$ of finite type satisfying the following:
\begin{enumerate}
\item For every finite generated $\Lambda$-module $M$, there is a graded ring homomorphism $$\phi_M: H \longrightarrow \Ext_{\Lambda}^{*}(M,M).$$
\item For each pair $(X,Y)$ of finite generated $\Lambda$-modules, the scalar actions from $H$ on $\Ext_{\Lambda}^{*}(X,Y)$ via $\phi_X$ and $\phi_Y$ coincide, and $\Ext_{\Lambda}^{*}(X,Y)$ is a  finite generated $\Lambda$-module.
\end{enumerate}
Bergh proves that $c(\Lambda) +1 \le \mathrm{repdim}\, \Lambda $ for a self-injective algebra satisfying assumption {\bf Fg} \cite{ber}.
The following corollary follows immediately.

\begin{cor}\label{repdim}
Let $\Lambda$ be a Koszul self-injective algebra and $\tilde{\Lambda}^{\sigma}= \Lambda\ltimes D\Lambda^{\sigma}$ be its twisted trivial extension.
If $\tilde{\Lambda}^{\sigma}$ satisfies the assumption {\bf Fg}, and if the complexity of $\Lambda$ is $r$, then $$\mathrm{repdim}\, \tilde{\Lambda}^{\sigma} \ge r + 2.$$
\end{cor}

%Let $ll(\Lambda)$ be  the Loewy length of $\Lambda$.
By Theorem \ref{main}, we have that if $\Lambda$ is a graded self-injective algebra, then the Loewy length of $\tilde{\Lambda}^{\sigma}$ is $1$ plus the Loewy length of $\Lambda$.
Apply Theorem 3.2 of \cite{ber}, one gets

\begin{cor}\label{repdima}
Let $\Lambda$ be a Koszul self-injective algebra and $\tilde{\Lambda}^{\sigma}= \Lambda\ltimes D\Lambda^{\sigma}$ be its twisted trivial extension.
Assume that both algebras satisfy the assumption {\bf Fg} of \cite{ber}
If  the Loewy length of $\Lambda$ is $c(\Lambda)+1$, then $$\mathrm{repdim}\, \tilde{\Lambda}.^{\sigma} = \mathrm{repdim}\, {\Lambda} +1.$$
\end{cor}

It is natural to conjecture that the representation dimension of a twisted trivial extension of a self-injective algebra increases by one in general.

\begin{exa}\label{exterior}
{\em An exterior algebra of $m$-dimensional vector space is  Koszul self-injective algebra of Loewy length $m+1$.
It has representation dimension $m+1$ \cite{ro} and complexity $m$.
Exterior algebra satisfies the assumption {\bf Fg}.
By Example \ref{skewcentral}, the exterior algebra of an $m+1$-dimensional vector is isomorphic to a twisted trivial extension of an exterior algebra of its $m$-dimensional subspace.
So Theorem \ref{complexity} and Corollary \ref{repdima} tell us how the complexities of the exterior algebras and the representation dimensions grow as the dimensions of the vector spaces grow.
}\end{exa}

\section{Extensions of Koszul Artin-Schelter Regular Algebras}

A graded algebra $\Gamma =  \Gamma_0+ \Gamma_1+ \Gamma_2+ \cdots$ is called an {\em Artin-Schelter regular algebra} of dimension $l$ if the following conditions are satisfied for the positive integer $l$:
\begin{enumerate}
 \item Each graded simple $\Gamma$-module has projective dimension $l$;
 \item For each graded simple $\Gamma$-module $S$ and for $0 \le t < l$, $\Ext_{\Gamma}^{t}(S, \Gamma) =0 $;

 \item $\{ \Ext_{\Gamma}^{m}(S, \Gamma) |S \mbox{ graded} \mbox{ simple } \Gamma-\mbox{module} \}$ is a complete set of graded simple $\Gamma^{op}$-modules.
\end{enumerate}
\noindent Such algebra is also called general Auslander regular algebra in \cite{gmt}.

For a graded algebra $\Gamma  =  \Gamma_0+ \Gamma_1+ \Gamma_2+ \cdots$, its Yoneda  algebra is the vector space $E(\Gamma)= \bigoplus_{i=1}^{\infty} \Ext_{\Gamma}^t( \Gamma_0 , \Gamma_0 )$ with the multiplication defined by the Yoneda product.
It is shown in \cite{sm, m2} that for a Koszul Artin-Schelter regular algebra $\Gamma$ of global dimension $l$, its  Yoneda algebra $\Lambda= E(\Gamma)$ is a Koszul self-injective algebra of Loewy length $l+1$, and verse versa.
The functors $F = \Ext^{\bullet}_{\Gamma} (\quad, \Gamma_0)$ and $G = \Ext^{\bullet}_{\Lambda} ( \quad, \Lambda_0) $ are dualities between the Koszul modules of Artin-Schelter regular algebras $\Gamma$ and the Koszul self-injective algebras $\Lambda$ \cite{gm}.
All the projective $\Lambda$-module have the same Loewy length, which is $1$ plus the projective dimension of simple $\Gamma$-modules.
So the Loewy length of $\Lambda$ is  $1$ plus the global dimension of $\Gamma$ \cite{m2}.
We also know that for a Koszul $\Gamma$-module $M$, the Gelfand-Kirilov dimension of $M$ is the same as the complexity of the $\Lambda $-module $F(M)$ \cite{g4}.

We have the following version of Theorem \ref{complexity} for Koszul Artin-Schelter regular algebras:
\begin{thm}\label{as}
Let $\Gamma$ be a Koszul Artin-Schelter regular algebra of global dimensional $l$ and Gelfand-Kirilov dimension $c$, then $E(E(\Gamma)\ltimes DE(\Gamma)^{\sigma})$ is a Koszul Artin-Schelter regular algebra of global dimensional $l+1$ and Gelfand-Kirilov dimension $c+1$.
\end{thm}

Note that $E(\Gamma)\simeq (\Gamma^!)^{op}$ is just the opposite algebra of the quadratic dual of $\Gamma$.
Assume that $\Gamma$ is given by the bound quiver $(Q , \theta)$, with relation set $\theta$, then $E(\Gamma)$ is given by the bound quiver $(Q , \rho)$, and we may take $\rho$ as a basis of the dual space of the subspace spanned by $\theta$ in the subspace of the path algebra $kQ$ spanned by the paths of length $2$.
If $(Q,\rho)$ is a stable bound quiver of Loewy length $l+1$ with Nakayama translation $\tau$, we also call $(Q, \theta)$ a stable bound quiver of Loewy length $l+1$ with Nakayama translation $\tau$.
We have the following version of Theorem \ref{main}, by taking the dual in Theorem \ref{main}, Proposition \ref{twistR} and Theorem \ref{tcomplexity}.

\begin{thm}\label{dualmain}
Let $\Gamma$ be a graded Artin-Schelter regular algebra given by a bound quiver $Q=(Q, \theta)$.
Let $\tilde{\Gamma}^{\sigma} =E(E(\Gamma)\ltimes DE(\Gamma)^{\sigma})$ for an automorphism $\sigma$ of $E(\Gamma)$, and let $\tilde{Q}^{\sigma}=(\tilde{Q},\tilde{\theta}^{\sigma})$ be a bound quiver of $\tilde{\Gamma}^{\sigma}$.
If $Q$ is a stable bound quiver of Loewy length $l+1$ with the  Nakayama translation $\tau=\tau_{E(\Gamma)}$. Then
\begin{enumerate}
\item $(\tilde{Q},\tilde{\theta}^{\sigma})$ is a stable translation quiver of Loewy length $l+2$ with trivial Nakayama translation, and $\tilde{Q}$ is obtained from $Q$ by adding an arrow $\alpha_i: i \longrightarrow \tau i$ for each vertex $i\in Q_0$;

\item The relation set $\tilde{\theta}^{\sigma}$ can be taken as $\tilde{\theta}^{\sigma}= \theta\cup \theta'$ with $$\theta'= \{\alpha_j \beta + \nu\sigma^{-1}(\beta)\alpha_i | \beta: i\longrightarrow j \in Q_1\}.$$

\item If the Gelfand-Kirilov dimension of $\Gamma$ is $c$, then the Gelfand-Kirilov dimension of $\tilde{\Gamma}^{\sigma}$ is $c+1 $.
\end{enumerate}\end{thm}

Now assume that $\Gamma$ is {\em connected} in the sense that $\Gamma_0 = k$.
Denote by $\Gamma^e = \Gamma \otimes \Gamma^{op}$ the enveloping algebra.
$\Gamma$ is called a graded {\em $l$-Calabi-Yau algebra} \cite{gin} if the following conditions hold:
\begin{enumerate}
\item As a $\Gamma^e$ module, $\Gamma$ has a finitely generated projective resolution of finite length;

\item $\Ext_{\Gamma^e}^{t}(\Gamma, \Gamma^e) =0$ for $ t \neq l$ and$ \Ext_{\Gamma^e}^{l}(\Gamma, \Gamma^e) \simeq A$ as $\Gamma^e$ modules.
\end{enumerate}
Calabi-Yau algebra is an important class of Artin-Schelter regular algebras, and found applications in many fiels of mathematics, such as Homological mirror symmetry \cite{Kon}, cluster algebras and cluster categories \cite{IR, KR}, string theory \cite{Laz}, conformal field theory \cite{Co}, and  Donaldson-Thomas invariants on Calabi-Yau threefold in geometrical invariant theory \cite{MR, Se}.

It is known that an Artin-Schelter regular algebra $\Gamma$ is $l$-Calabi-Yau if and only if $E(\Gamma)$ is graded symmetric(see Proposition 2.10 of \cite{wz}).

Now assume that $\Gamma$ is connected.
Write $ \tilde{ \Gamma }^{ \sigma } = E ( E ( \Gamma ) \ltimes D E (\Gamma)^{ \sigma } )$.
The quiver of a connected graded algebra is a quiver with only one vertex and arrows are just loops.
So Theorem \ref{dualmain} tells us that there is just one more loop $y$ in the quiver $\tilde{Q}$ of $ \tilde{ \Gamma }^{ \sigma }$.

Since $\Gamma$ is Koszul of global dimension $l$, $E(\Gamma)$ has Loewy length $l+1$.
Under the automorphism $\varepsilon'$ on $E(\Gamma)$, $ E(\Gamma) \ltimes D E( \Gamma )^{ \varepsilon' }$ is graded symmetric by Example \ref{gradeds}.
The following Proposition follows from Proposition 2.10 of \cite{wz}.
\begin{prop}\label{CY}
If $ \Gamma $ is a connected Koszul Artin-Schelter regular algebra of dimension $l$, then $ \tilde{ \Gamma } = E ( E ( \Gamma ) \ltimes D E ( \Gamma )^{\varepsilon'} )$ is a Koszul $l+1$-Calabi-Yau algebra.
\end{prop}

For $\sigma_0 = \varepsilon \nu $, the Yoneda algebra of $ \tilde{ \Gamma }^{ \sigma_0 }$ is described in Example \ref{skewcentral}.
So the new generator $y$ satisfies the condition $$ xy=yx $$ for any $x\in \Gamma$.
Thus $\tilde{\Gamma}^{\sigma_0}$ is just the central extension of $\Gamma$.
With $\Gamma[\sigma]= \tilde{\Gamma}^{\sigma\sigma_0}$, we have the following theorem.
\begin{thm}\label{extfamily}
Let $\Gamma$ be a Koszul Artin-Schelter regular algebra of global dimensional $l$ and Gelfand-Kirilov dimension $c$, then we have a family $\{\Gamma[\sigma] | \sigma \in \mathrm{Aut}\, E(\Gamma)\}$ of Koszul Artin-Schelter regular algebras of global dimensional $l+1$ and Gelfand-Kirilov dimension $c+1$ parameterized by the group $Aut_g(E(\Gamma))$ of the graded automorphisms of the Yoneda algebra of $\Gamma$.

If $\Gamma$ is connected then $\Gamma[\nu^{-1}\varepsilon^l]$ is $l+1$-Calabi-Yau and $\Gamma[1]$ is a central extension of $\Gamma$.
\end{thm}

We call each algebra $\Gamma[\sigma]$ a {\em  generalized twisted central extension} of $\Gamma$ with the twister $\sigma$, and call the Calabi-Yau one a {\em Calabi-Yau extension}.

{}


\begin{thebibliography}{}


\bibitem{AS} Artin, M. and Schelter, W. F.: Graded Algebras of Global Dimension $3$, Adv. Math. 66 (1987), 171--216.


\bibitem{ATV1} Artin, M. Tate, J. and Van den Bergh, M. Some algebras associated to automorphisms of elliptic curves,  {\it The Grothendieck Festschrift}, 33--85, Birkh\"{a}user, Boston, 1990.

\bibitem{ATV2} Artin, M. Tate, J. and Van den Bergh, M.: Modules over regular algebras of dimension $3$, Invent. Math., 106 (1991) 335--388.

\bibitem{ars}  Auslander,  M.,  Reiten,  I. and  Smal\o , S.: {\em Representation theory of artin algebras}, Cambridge Studies in Advanced Math. 36 (1995), Cambridge Univ. Press.

\bibitem{bgs} Beilinson, A. Ginsberg, V. and Soergel, W. : Koszul duality patterns in Representation theory, J. of AMS 9 no.2 (1996) 473-527

\bibitem{ber} Bergh,P.A.: Representation dimension and finitely generated cohomology, Adv. Math. 219 (2008), no. 1, 389-400.

\bibitem{CS} Cassidy, T. Shelton,  B.:   PBW-deformation theory and regular central extensions,  J. Reine Angew. Math., 610 (2007), 1--12.

\bibitem{Co} Costello, K.: Topological conformal field theories and Calabi-Yau categories, Adv. Math. 210 (2007),  165--214.

\bibitem{dt} de la Pe\~na J. A. and Takane M.: On the Coxeter transformations of a wild algebra, Archiv der Math. 63(1994) 128 - 135.

\bibitem{gls} Ch.Geiss, B. Leclerc, and J. Schr\"oer: Rigid modules over preprojective algebras, Invent. Math. vol 165(2006), no.3, 589-632

\bibitem{gin} Ginzburg, V.: Calabi-Yau algebras, preprint arXiv:0612139

\bibitem{gm} Green, E. L. and Mart\'inez-Villa, R. Koszul and Yoneda algebras, Representation theory of algebras, CMS Conference Proceedings vol 19(1996) 247-298.

\bibitem{g4}  Guo, J. Y. and Wu Q.: Loewy matrix, Koszul cone and applications,  Comm. in algebra, 28(2000) no.2 925-941


\bibitem{gum}  Guo, J. Y., Mart\'inez-Villa: Algebra pairs associated to McKay quivers. Comm. in Algebra 30 (2002), no. 2, 1017--1032.

\bibitem{gmt}  Guo, J. Y., Mart\'inez-Villa and Takane M.: Koszul generalized Auslander regular algebras. in Algebras and Modules II, Canadian Mathematical Society Conference Proceedings 24(1998) 263-283


\bibitem{gcv}  Guo, J. Y.: On  McKay quivers and covering spaces. arXiv:1002.1768

\bibitem{g10}  Guo, J. Y.: Coverings and Truncations of Graded Self-injective Algebras, preprint arXiv:1002.4910v2

\bibitem{iy} Iyama, O.:	Cluster tilting for higher Auslander algebras. 	 arXiv:0809.4897

\bibitem{IR} Iyama, O. and Reiten, I. :  Fomin-Zelevinsky mutation and tilting modules over Calabi-Yau algebras, Amer. J. Math. 130 (2008), 1087--1149.

\bibitem{KR} Keller, B. and Reiten, I.: Cluster-tilted algebras are Gorenstein and stably Calabi-Yau, Adv. Math. 211 (2007), 123--151.

\bibitem{Kon} M. Kontsevich, M.: Homological algebra of mirror symmetry,  {\it Proc. Internat. Congress of Mathematicians}, (Z\"{u}rich, 1994), Birkh\"{a}user, Basel, 1995,  120--139.

\bibitem{Laz} Lazaroiu, C. I. : String field theory and brane superpotentials,  J. High Energy Phys. (2001), Paper 18, 40 pp.

\bibitem{lsv} Le Bruyn, L.,Smith, S.P., Van den Bergh, M.: Central extensions of three-dimensional Artin-Schelter regular algebras. Math. Z. 222, 171¨C212 (1996)


\bibitem{Lev} Levasseur, T.:  Some properties of noncommutative regular graded rings. Glasgow Math. J. 34 (1992),  277--300.

\bibitem{LPWZ} Lu, D.M. Palmieri, J.H. Wu, Q.-S. and Zhang, J. J. : Regular algebras of dimension 4 and their $A_\infty $-Ext-algebras, Duke Math. J. 137 (2007) 537-584.

\bibitem{m2} Mart\'{\i}nez-Villa, R.: Graded, selfinjective, and Koszul algebras. J. Algebra 215 1 (1999), pp. 34¨C72.

\bibitem{MR} Mozgovoy S.: and Reineke, M.: On the noncommutative Donaldson-Thomas invariants arising from brane tilings, Adv. Math. 223 (2010), 1521--1544.

\bibitem{ro} Rouquier, R.: Representation dimension of exterior algebras, Invent. Math. 165 (2006), 357-367.

\bibitem{sm} Smith, P.: Some finite dimensional algebras related to elliptic curves, Representation theory of algebras and related topics,   CMS Conference Proceedings vol 19 315-348.

\bibitem{SSt} Smith, S.P. and Stafford, J.T. : Regularity of the four dimensional Sklyanin algebra, Compositio Math. 83 (1992) 259--289.

\bibitem{Se} Segal, E. : The A$_\infty$ deformation theory of a point and the derived categories of local Calabi-Yaus, J. Algebra 320 (2008), 3232--3268

\bibitem{wz} Q.S wu and C. Zhu.: Poincar\'e-Birkhoff-Witt Deformation of Koszul Calabi-Yau algebras. Algebras and Representation Theory (to appear)
\end{thebibliography}
\end{document}